\documentclass[a4paper,10pt,leqno]{article}
\pdfoutput=1

\usepackage[english]{babel}
\usepackage{lipsum}

\newcommand{\name}[1]{{#1}}
\newcommand{\ie}{i.e.,~}
\newcommand{\eg}{e.g.,~}
\newcommand{\svp}{\textsc{Svp}}
\newcommand{\cldl}{Cl-\textsc{Dlp}}
\newcommand{\approxpb}[1]{Approx-{#1}}

\usepackage{enumerate}
\usepackage{enumitem}
\setlist{%
  topsep=.2\baselineskip,%
  parsep=0pt,%
  itemsep=0\baselineskip,%
}
\usepackage[dvipsnames]{xcolor}

\usepackage{array}
\usepackage{multirow}
\newcolumntype{C}{>{$}c<{$}}
\newcolumntype{L}{>{$}l<{$}}
\newcolumntype{R}{>{$}r<{$}}
\usepackage[format=hang,labelfont=sc,labelsep=endash]{caption}

\usepackage{pgf,pgfcore,pgfkeys,tikz}
\usetikzlibrary{fit,matrix,positioning}

\usepackage{amsmath,amsthm}
\usepackage{amssymb,latexsym}
\usepackage{stmaryrd}

\newcommand\Z{\mathbb{Z}}
\newcommand\Q{\mathbb{Q}}
\newcommand\F{\mathbb{F}}
\newcommand{\range}[2]{\llbracket {#1}, {#2} \rrbracket}
\newcommand{\nvars}[3][1]{#2_{#1}, \dotsc, #2_{#3}} % \nvars[j]{x}{n}: x_j, ..., x_n

\newcommand{\lnorm}[2][\empty]{\lVert{#2}\rVert_{#1}}

\newcommand{\NormOp}{\mathcal{N}}
\newcommand{\norm}[1]{\NormOp(#1)}
\newcommand{\normrel}[2]{\NormOp_{#1/#2}}
\DeclareMathOperator{\TraceOp}{Tr}

\newcommand{\QQ}[1]{K_{#1}}
\newcommand{\cf}[1]{\QQ{#1}} % Def. ob
\newcommand{\z}[1]{\zeta_{#1}}
\newcommand{\realfield}[1]{#1^{+}}

\newcommand{\classnb}[1]{h_{#1}}
\newcommand{\classnbreal}[1]{\classnb{#1}^{+}}
\newcommand{\classnbminus}[1]{\classnb{#1}^{-}}
\newcommand{\maxorder}[1]{\mathcal{O}_{#1}}
\newcommand{\ideal}[1]{\mathfrak{#1}}
\newcommand{\pideal}[1]{\langle #1 \rangle}
\newcommand{\units}[1]{\maxorder{#1}^{\times}}

\newcommand{\circunits}[1]{C_{#1}}
\newcommand{\ufund}[1]{\varepsilon_{#1}}
\DeclareMathOperator{\Log}{Log}
\newcommand{\ulog}{\Log_{\infty}}

\DeclareMathOperator{\Gal}{Gal}
\newcommand{\GG}[1]{G_{#1}}
\newcommand{\si}[2]{\sigma_{#2,#1}}
\newcommand{\res}{\operatorname{res}}
\newcommand{\cor}{\operatorname{cor}}
\newcommand{\cplxconj}[1]{\tau}

\newcommand{\St}[1]{\mathcal{S}_{#1}}
\newcommand{\sti}[1]{\St{#1}} % Def. ob
\newcommand{\Stp}[1]{\mathcal{S}'_{#1}}
\newcommand{\Stpp}[1]{\mathcal{S}''_{#1}}
\newcommand{\R}{\mathcal{R}}
\newcommand{\A}{\mathcal{A}}

\newcommand{\LL}{\ideal{L}}
\newcommand{\Su}{\mathcal{S}}

\usepackage{chngcntr}

\usepackage{hyperref}

\hypersetup{
        final,                    % Turn on hyperref even in draft version          [default:cls option]
        hypertexnames = false,    % Workaround "already defined identifiers" warns  [default:true]
        bookmarksnumbered,        % Numbering appears in bookmarks names            [default:false]
        colorlinks = true,
        citecolor = RoyalBlue,
        linkcolor = RoyalPurple,
        urlcolor  = ForestGreen,
}

\usepackage{pgf,pgfcore,pgfkeys}

\usepackage[capitalize,nosort]{cleveref}

\makeatletter

\newcommand{\setcref@def@abbrevspace}{\,}
\newcommand{\setcref@def@fullspace}{~}

\newcommand{\setcref@def@rangeconjunction}{--}
\newcommand{\setcref@def@pairconjunction}{ and}
\newcommand{\setcref@def@middleconjunction}{,}
\newcommand{\setcref@def@lastconjunction}{, and}

\pgfkeys{
    /@setcref/.is family, /@setcref,
    default/.style = {%
        abbrevspace = \setcref@def@abbrevspace,
        fullspace   = \setcref@def@fullspace,
        delim       = none,
        forcelow    = false,
    },
    abbrevspace/.store in  = \setcref@abbrevspace,
    fullspace/.store in    = \setcref@fullspace,
    delim/.style           = {delims/#1},
    forcelow/.style        = {ifforcelow/#1},
    /@setcref/delims/.cd,
        none/.code     = {\def\setcref@delim@begin{}%
                          \def\setcref@delim@end{}},
        parent/.code   = {\def\setcref@delim@begin{(}%
                          \def\setcref@delim@end{)}},
        brack/.code    = {\def\setcref@delim@begin{$\{$}%
                          \def\setcref@delim@end{$\}$}},
        sqbrack/.code  = {\def\setcref@delim@begin{[}%
                          \def\setcref@delim@end{]}},
    /@setcref/ifforcelow/.cd,
        true/.code     = {\edef\setcref@capitalize{}},
        false/.code    = {\edef\setcref@capitalize{\if@cref@capitalise\noexpand\MakeUppercase\fi}},
}

\newcommand{\setcref}[6][]{%
        \pgfkeys{/@setcref, default, #1}%

        \expandafter\edef\csname setcref@#2@capitalize\endcsname%
                         {\unexpanded\expandafter{\setcref@capitalize}}
        \expandafter\edef\csname setcref@#2@abbrevspace\endcsname%
                         {\unexpanded\expandafter{\setcref@abbrevspace}}
        \expandafter\edef\csname setcref@#2@fullspace\endcsname%
                         {\unexpanded\expandafter{\setcref@fullspace}}
        \expandafter\edef\csname setcref@#2@delim@begin\endcsname%
                         {\unexpanded\expandafter{\setcref@delim@begin}}
        \expandafter\edef\csname setcref@#2@delim@end\endcsname%
                         {\unexpanded\expandafter{\setcref@delim@end}}

        \crefname{#2}{%
            \csname setcref@#2@capitalize\endcsname#5%
          }{%
            \csname setcref@#2@capitalize\endcsname#6%
        }
        \Crefname{#2}{#3}{#4}

        \crefformat{#2}{%
            \csname cref@#2@name\endcsname%
            \csname setcref@#2@abbrevspace\endcsname\csname setcref@#2@delim@begin\endcsname%
            ##2##1##3%
            \csname setcref@#2@delim@end\endcsname%
        }
        \Crefformat{#2}{%
            \csname Cref@#2@name\endcsname%
            \csname setcref@#2@fullspace\endcsname\csname setcref@#2@delim@begin\endcsname%
            ##2##1##3%
            \csname setcref@#2@delim@end\endcsname%
        }

        \crefrangeformat{#2}{%
            \csname cref@#2@name@plural\endcsname%
            \csname setcref@#2@abbrevspace\endcsname\csname setcref@#2@delim@begin\endcsname%
            ##3##1##4%
            \setcref@def@rangeconjunction% [opt?]
            ##5\crefstripprefix{##1}{##2}##6%
            \csname setcref@#2@delim@end\endcsname%
        }
        \Crefrangeformat{#2}{%
            \csname Cref@#2@name@plural\endcsname%
            \csname setcref@#2@fullspace\endcsname\csname setcref@#2@delim@begin\endcsname%
            ##3##1##4%
            \setcref@def@rangeconjunction% [opt?]
            ##5\crefstripprefix{##1}{##2}##6%
            \csname setcref@#2@delim@end\endcsname%
        }
        
        \crefmultiformat{#2}{% First
            \csname cref@#2@name@plural\endcsname%
            \csname setcref@#2@abbrevspace\endcsname\csname setcref@#2@delim@begin\endcsname%
            ##2##1##3%
            \csname setcref@#2@delim@end\endcsname%
          }{% Second
            \setcref@def@pairconjunction%
            \csname setcref@#2@abbrevspace\endcsname\csname setcref@#2@delim@begin\endcsname%
            ##2##1##3%
            \csname setcref@#2@delim@end\endcsname%
          }{% Middle
            \setcref@def@middleconjunction%
            \csname setcref@#2@fullspace\endcsname\csname setcref@#2@delim@begin\endcsname%
            ##2##1##3%
            \csname setcref@#2@delim@end\endcsname%
          }{% Last
            \setcref@def@lastconjunction%
            \csname setcref@#2@abbrevspace\endcsname\csname setcref@#2@delim@begin\endcsname%
            ##2##1##3%
            \csname setcref@#2@delim@end\endcsname%
        }
        \Crefmultiformat{#2}{% First
            \csname Cref@#2@name@plural\endcsname%
            \csname setcref@#2@fullspace\endcsname\csname setcref@#2@delim@begin\endcsname%
            ##2##1##3%
            \csname setcref@#2@delim@end\endcsname%
          }{% Second
            \setcref@def@pairconjunction%
            \csname setcref@#2@abbrevspace\endcsname\csname setcref@#2@delim@begin\endcsname%
            ##2##1##3%
            \csname setcref@#2@delim@end\endcsname%
          }{% Middle
            \setcref@def@middleconjunction%
            \csname setcref@#2@fullspace\endcsname\csname setcref@#2@delim@begin\endcsname%
            ##2##1##3%
            \csname setcref@#2@delim@end\endcsname%
          }{% Last
            \setcref@def@lastconjunction%
            \csname setcref@#2@abbrevspace\endcsname\csname setcref@#2@delim@begin\endcsname%
            ##2##1##3%
            \csname setcref@#2@delim@end\endcsname%
        }
}

\newcommand{\crefext}[2]{\csname cref@#1@format\endcsname{#2}{}{}}
\newcommand{\crefrangeext}[3]{\csname crefrange@#1@format\endcsname{#2}{#3}{}{}{}{}}

\makeatother

\setcref[abbrevspace=\empty]{section}    {Section}{Sections}{\S}{\S\S}
\setcref[forcelow=true,%                 
         abbrevspace=\empty]{page}       {Page}{Pages}{p.}{pp.}
\setcref                    {table}      {Table}{Tables}{tab.}{tab.}
\setcref[delim=parent,%
        abbrevspace=\,]     {equation}   {Equation}{Equations}{\!\!}{\!\!}
\setcref[delim=parent,%
        abbrevspace=\,]     {eqext}      {Equation}{Equations}{\!\!}{\!\!}
\setcref                    {rem}        {Remark}{Remarks}{rem.}{rems.}
\setcref                    {exercice}   {Exercice}{Exercices}{ex.}{ex.}
\setcref                    {defin}      {Definition}{Definition}{def.}{def.}
\setcref                    {lem}        {Lemma}{Lemmas}{lem.}{lem.}
\setcref                    {prop}       {Proposition}{Propositions}{pr.}{pr.}
\setcref                    {thm}        {Theorem}{Theorems}{th.}{th.}
\setcref                    {corol}      {Corollary}{Corollaries}{cor.}{cor.}
\setcref                    {algorithm}  {Algorithm}{Algorithms}{alg.}{alg.}
\newtheorem{thm}  {Theorem}
\newtheorem{corol}[thm]{Corollary}
\newtheorem{prop} [thm]{Proposition}
\newtheorem{lem}  [thm]{Lemma}

\theoremstyle{definition}

\newtheorem*{xrem}{Remark}

\counterwithin{thm}{section}
\counterwithin{table}{section}

\newcommand{\prometheus}{The first-named author is supported by the European Union PROMETHEUS project (Horizon 2020 Research and Innovation Program, grant 780701).}
\newcommand{\mail}[1]{\href{mailto:#1}{\texttt{#1}}}

\title{\Large \bf A short basis of the Stickelberger ideal\\ of a cyclotomic field%
\thanks{\prometheus}}
\author{\normalsize Olivier Bernard\textsuperscript{1,2} and Radan Ku\v{c}era\textsuperscript{3}\\[7pt]
  \small \textsuperscript{1} Univ Rennes, CNRS, IRISA, France\\[-3pt]
  \small \mail{olivier.bernard@irisa.fr}\\[-1pt]
  \small \textsuperscript{2} Thales, Laboratoire Chiffre, Gennevilliers, France\\[-1pt]
  \small \textsuperscript{3} Faculty of Science, Masaryk University, Brno, Czech Republic\\[-3pt]
  \small \mail{kucera@math.muni.cz}\\[7pt]
}

\date{\normalsize\today}
\newcommand{\lkw}{Cyclotomic fields, Stickelberger ideal, short basis, relative class number}
\providecommand{\keywords}{\textbf{Keywords:} }

\hypersetup{
  pdftitle    = {A short basis of the Stickelberger ideal of a cyclotomic field},
  pdfauthor   = {Olivier Bernard and Radan Kucera}, % No special character in pdf metadata... :(
  pdfsubject  = {A short basis of the Stickelberger ideal of a cyclotomic field, ArXiv Pre-Print [math.NT]}, % ...
  pdfkeywords = {\lkw}, % ...
  pdflang     = {English},
}

\begin{document}
\maketitle
\begin{abstract}
    We exhibit an explicit \emph{short} basis of the Stickelberger ideal of cyclotomic fields of any conductor $m$, \ie a basis containing only short elements. By definition, an element of $\Z[\GG{m}]$, where $\GG{m}$ denotes the Galois group of the field, is called short whenever it writes as $\sum_{\sigma\in\GG{m}} \varepsilon_{\sigma}\sigma$ with all~$\varepsilon_\sigma\in\{0,1\}$.
    One ingredient for building such a basis consists in picking wisely generators~$\alpha_m(b)$ in a large family of short elements.
    
    As a direct practical consequence, we deduce from this short basis an \emph{explicit} upper bound on the relative class number, that is valid for \emph{any} conductor.
    This basis also has several concrete applications, in particular for the cryptanalysis of the Shortest Vector Problem on Ideal lattices. % [Grammar] If has is the only verb (no past / auxiliary / etc) then also is before.

    \bigskip\noindent\keywords\lkw\bigskip
\end{abstract}

% ----------------------------------------------------------------------------------------------

% Context: cryptanalysis of id-SVP
% --------------------------------------
% 1. Avent of quantum computers => new problems => SVP in ideal lattices, cyclotomic fields.

The eventuality of achieving to build a large scale quantum computer in the next few decades has driven the cryptographic community to consider new mathematical problems upon which to base so-called \emph{post-quantum} cryptosystems.
Among many proposals, some of the most promising ultimately base their security on the hardness of the Approximate Shortest Vector Problem (\approxpb{\svp}) in algebraically structured Euclidean lattices, that offers a nice trade-off between security and efficiency.
In fact, a popular choice is to consider fractional ideals in some cyclotomic field~$\cf{m}=\Q[\z{m}]$ of conductor~$m\not\equiv 2\!\pmod 4$, \eg$m=2048$. Such an ideal~$\ideal{b}$ can be viewed as a Euclidean lattice under the Minkowski embedding, and the Approximate Ideal-\svp{} consists in finding~$x\in \ideal{b}$ such that the induced Euclidean norm~$\lnorm[2]{x}$ is close to the smallest possible one.

% 2. Cryptanalytic effort:
% Remark that in the prime-power case, good basis of Log U.
% Extend to all ideals using the Stickelberger lattice

% Cryptanalytic effort to take advantage of this additional number field structure for breaking ideal SVP.
In the last decade, there has been a significant cryptanalytic effort trying to benefit from this additional algebraic structure to solve \approxpb{Ideal-\svp}, giving rise to a long series of works (\cite{CGS14,CDPR16,CDW17,DPW19,PHS19,BR20,CDW21}).
All start from a solution to the Class Group Discrete Logarithm Problem (\cldl), which is, given a fixed set of finite places corresponding to prime ideals~$\bigl\{\nvars{\ideal{p}}{k}\bigr\}$ of~$\cf{m}$, and any challenge ideal~$\ideal{b}$,%
\footnote{Actually, for the problem to have a solution, the challenge ideal shall be chosen such that its class in the class group of~$\cf{m}$ belongs to the subgroup generated by the classes of the~$\ideal{p}_i$'s.}\ %
to find~$\alpha\in \cf{m}$ and~$\nvars{e}{k}\in\Z$ such that
\begin{equation*}
    \pideal{\alpha} = \ideal{b} \cdot \prod_{1\leq i\leq k} \ideal{p}_i^{e_i}.
\end{equation*}
In a quantum world, it appears this problem is not hard to solve (\cite{EHKS14,BS16}), so the most difficult part of these cryptanalyses resides in reducing the Euclidean norm of~$\alpha$.
When the challenge~$\ideal{b}=\pideal{\alpha}$ is guaranted to be principal, so that it is possible to consider an empty set of finite places, the traditional method consists in using the log-unit lattice of~$\cf{m}$, \ie the image of the unit group~$\units{\cf{m}}$ under the logarithmic embedding~$\ulog$ of~$\cf{m}$. 
Indeed, 
if the closest vector to~$\ulog \alpha$ in~$\ulog \units{\cf{m}}$ is~$\ulog \ufund{}$,~$\ufund{}\in\units{\cf{m}}$, then~$\alpha/\ufund{}$ is hopefully the smallest generator of~$\ideal{b}$.
This principle can be extended to general ideals \cite{PHS19,BR20} by using the log-$\Su$-unit lattice under the~$\Su$-logarithmic embedding \cite[\crefext{section}{3},~\crefext{page}{98}]{Nar04}, where~$\Su$ contains the chosen finite places as well as all infinite places.

In general lattices, finding the closest vector to any target is a well-known hard problem. In the above case though, it has been noticed in \cite{CGS14} and thereafter rigorously proven \cite{CDPR16,CDW21} that logarithmic embeddings from the set~$\circunits{m}$ of circular units \cite[\crefext{section}{8}]{Was97} yield a sufficiently good basis of a sublattice of~$\ulog \units{\cf{m}}$ of relatively small finite index. A key property of these vectors is that they are small compared to the regulator of~$\cf{m}$. An explicit set of independent generators of~$\circunits{m}$ has been given for any conductor~$m$ in \cite[\crefext{thm}{2}]{GK89} and independently in \cite[\crefext{thm}{6.1}]{K92}. 
No such result for the quality of some explicit basis has been proven for log-$\Su$-unit lattices for a non-empty set of finite places, though limited experimental evidence in the prime conductor case \cite{BR20} tend to show that this phenomenon still holds.

Furthermore, by Stickelberger's theorem, the Stickelberger ideal~$\sti{m}$ of~$\cf{m}$ annihilates its class group, so it was proposed in \cite{CDW17,CDW21} to use these free relations to help reducing the algebraic norm of the \cldl{} solution.
More precisely, since by~\cite{S78}~$(1-\cplxconj{m})\sti{m}$, viewed as a $\Z$-module, has full rank in~$(1-\cplxconj{m})\Z\bigl[\GG{m}\bigr]$, where~$\GG{m}=\Gal(\cf{m}/\Q)$ and $\cplxconj{m}\in\GG{m}$ is induced by complex conjugation, it is a lattice of class relations for the relative class group. %
Therefore, choosing a challenge ideal $\ideal{b}$ and 
prime ideals for the \cldl{} in the relative class group, \eg exactly one Galois orbit~$\bigl\{\ideal{p}^{\sigma}\bigr\}$ for all~$\sigma \in \GG{m}$, it is once again possible to express the reduction of a solution~$\pideal{\alpha}=\ideal{b}\cdot\ideal{p}^{\sum e_{\sigma}\sigma}$ as a closest vector problem in~$(1-\cplxconj{m})\sti{m}$, where the target is the vector~$(e_{\sigma}-e_{\tau\sigma})_{\sigma}$.
As noticed in \cite[\crefext{lem}{4.4}~and\,4.6]{CDW21}, this lattice contains many short elements, which \textit{in fine} yield a good description for finding sufficiently close vectors.
Note also that the plus part of the class group seems to be much smaller than the relative part,%
\footnote{This is backed up by several theoretical and computational observations, see \eg Weber's conjecture~$\classnbreal{2^e}=1$, Buhler, Pomerance and Robertson's conjecture for odd prime powers \cite{BPR04}, and Schoof's extensive calculations in~\cite[Tab.,\,\crefext{section}{4}]{Was97} and~\cite{Sch03}.}
hence every challenge $\ideal{b}$ can be reduced to this case by randomly searching for a small norm ideal $\ideal{c}$ such that the class of $\ideal{cb}$ belongs to the relative class group \cite[\crefext{algorithm}{5}]{CDW21}.

\paragraph{In praise of short Stickelberger bases.}
Unfortunately, while in the prime conductor case the exhibited set of short elements from \cite[\crefext{subsection}{4.2}]{CDW21} form a~$\Z$-basis
of~$\sti{m}$, in the general case this family is only known to generate~$\sti{m}$ as a~$\Z$-module.
This comes at the expense of constructing a linearly independent subset of vectors \cite[\crefext{corol}{2.2}]{CDW21}
that will only generate some full-rank sublattice, and should finally yield inferior approximation factors.
Worse, it is not even clear whether it is always possible to extract a basis from such a generating set, which may be crucial for some applications.

Another very important point is that the proof that the Stickelberger ideal annihilates the class group is completely explicit \cite[\crefext{section}{6.2}]{Was97}. Namely, for any prime ideal~$\ideal{p}$, and any~$\alpha\in\sti{m}$, it builds an explicit~$\gamma \in \cf{m}$ such that~$\pideal{\gamma} = \ideal{p}^{\alpha}$. However, if~$\alpha$ has 
even moderately large coefficients, this has an exponential impact on the height of~$\gamma$, that renders its computation rapidly intractable. On the contrary, having only short elements in the basis keeps the algebraic norm of the generators as low as possible, namely~$\norm{\ideal{p}}^{\varphi(m)/2}$. 
Explicitly computing Stickelberger generators is useful in at least two situations:
\begin{enumerate}
    \item the first one is when reducing the algebraic norm of the \cldl{} solution as in \cite{CDW21}, as knowing explicit generators prevents to perform a quantum step -- or, a classically costly step -- to recover the generator of the reduced ideal (see \cite[\crefext{thm}{5.1}]{CDW21} for the complete workflow);
    \item the second one occurs when one wants to use the knowledge of the Stickelberger relations to approach some log-$\Su$-unit lattice. Indeed, suppose the finite places of~$\Su$ correspond to one split Galois orbit~$\bigl\{\ideal{p}^{\sigma}\bigr\}$ for all~$\sigma \in \GG{m}$.
    Then, from a maximal set of independent real $\Su^+$-units, where the finite places of $\Su^+$ correspond to all relative norm ideals $\normrel{\cf{m}}{\realfield{\cf{m}}}\bigl(\ideal{p}^{\sigma}\bigr)$, adding explicit generators corresponding to a basis of the Stickelberger ideal, besides the absolute norm, yields a maximal set of independent $\Su$-units, at the much smaller cost of finding generators in the maximal real subfield. 
\end{enumerate}
In the latter case, note that knowing merely a short generating set of $\sti{m}$ instead of a $\Z$-basis is not sufficient %
to provide a full-rank family of independent $\Su$-units. Building a basis from such a generating set using the Hermite Normal Form would increase dramatically the height of the generators. Hence, having in the first place an \emph{explicit short basis} of $\sti{m}$ as a $\Z$-module is crucial here.

\paragraph{Historical results.}
The first explicitly known basis of~$\sti{m}$, viewed as a~$\Z$-module and for \emph{any} conductor~$m$, was given in \cite[\crefext{thm}{6.2}]{K92}, but elements of this basis have rather large coefficients.
In the prime conductor case, a short basis can be found in \cite[\crefext{thm}{9.3(i)}]{Sch10}, the shortness being proven in \cite[\crefext{exercice}{9.3}]{Sch10}.
This result has been extended to prime-power conductors in \cite{CDW17}, at the price of allowing slightly larger coefficients \cite[\crefext{lem}{4(2)}]{CDW17}. Finally, a large set of short \emph{generators} has been given in \cite[\crefext{section}{4.2}]{CDW21} in the general case for any conductor.

\paragraph{Contributions.}
In this work, our main result (see \cref{Short_basis}) is to provide the first explicit \emph{basis} of the Stickelberger ideal~$\sti{m}$ for any conductor~$m$, viewed as a~$\Z$-module, that is constituted \emph{only} of short elements, \ie elements of the form
\begin{equation*}
    \sum_{\sigma\in\GG{m}} a_{\sigma} \sigma \in \sti{m} \subset \Z\bigl[\GG{m}\bigr],\qquad\text{where $a_{\sigma}\in\{0,1\}$ for all~$\sigma$.}
\end{equation*}
Actually, besides the absolute norm element, all other members of this short basis have exactly~$\varphi(m)/2$ non-zero coordinates. In the prime conductor case, our short basis coincides with the basis given in \cite[\crefext{thm}{9.3(i)}]{Sch10}.
One ingredient of independent interest in the proof is \cref{Short}, which describes a large family of short elements of~$\sti{m}$ that encompasses the set from \cite[\crefext{section}{4.2}]{CDW21}, using a very simple arithmetic criterion in the spirit of \cite[\crefext{lem}{16.3}]{Was97} when $m$ is an odd prime power.
Picking wisely some elements~$\alpha_m(b)$ in this large family yields our proposed short basis.

We also show how to explicitly compute algebraic integers generating~$\mathfrak{L}^{\alpha_m(b)}$, for any unramified prime ideal~$\mathfrak{L}$ and any element~$\alpha_m(b)$ of our short basis. These generators can be expressed as Jacobi sums that turn out to be drastically more efficient to compute than the generators given \eg in \cite[\crefext{section}{6.2}]{Was97}.

Finally, a nice theoretical consequence of our result is to derive an explicit upper bound on the relative part~$\classnbminus{\cf{m}}$ of the class number of~$\cf{m}$. More precisely, for any conductor~$m \not\equiv 2\pmod 4$, \cref{Upper_bound} gives that
\begin{equation*}
    \classnbminus{\cf{m}} \leq 2^{1-a}\cdot \Bigl(\frac{\varphi(m)}{8}\Bigr)^{\varphi(m)/4},
\end{equation*}
where
\begin{equation*}
    a = \begin{cases}
        0         & \text{if~$m$ is a prime-power,}\\
        2^{t-2}-1 & \text{if~$m$ has~$t>1$ prime divisors.} 
    \end{cases}
\end{equation*}
To our knowledge, the best explicit upper bound on the relative class number which is valid for {any} conductor is given by \cite[\crefext{eqext}{6}]{Lou14}. However, whereas our bound is given by a simple formula and easy to manipulate, \name{Louboutin}'s bound is difficult to instantiate for comparison in the general case. As an example, the special case~$m = 4p$, where~$p\geq 3$ is an odd prime, is concretely treated in \cite[\crefext{thm}{2}]{Lou14}, which results in the following upper bound:
\begin{equation*}
    \classnbminus{\cf{4p}} \leq 8\sqrt{p}\cdot\Bigl(\frac{p}{16}\Bigr)^{(p-1)/2}.
\end{equation*}
We stress that in this example, this upper bound is sharper than ours.

We should also mention that the proof of our bound indirectly gives an algorithm to compute the relative class number by computing the determinant of some scaled Hadamard matrix: incidentally, this method seems to be significantly more efficient than when using the traditional analytic formula \cite[\crefext{thm}{4.17}]{Was97}, when the number~$t$ of prime factors of~$m$ is small.

\section{Notations and preliminaries}

For any integers $i, j$ with $i\leq j$, let $\range{i}{j}$ denote the set $\{k \in \Z;\ i\leq k\leq j\}$.
For any positive integer $m$,
let $\zeta_m=e^{2\pi i/m}$, let $\QQ{m}=\Q(\zeta_m)$ be the $m$th cyclotomic field and let
$\GG{m}=\Gal(\QQ{m}/\Q)$ be its Galois group.
Note that if $m$ is odd, we have $\QQ{2m}=\QQ{m}$ and $\GG{2m}=\GG{m}$.
For any $a\in\Z$, let
\begin{equation*}
\theta_m(a)=\sum_{\substack{0<s\le m\\(s,m)=1}}\left\langle-\frac{as}m\right\rangle\si{s}{m}^{-1}\in\Q[\GG{m}],
\end{equation*}
where $\langle x\rangle$ is the fractional part of a rational number $x$ (i.e. verifying~$0\le\langle x\rangle<1$ and~$x-\langle x\rangle\in\Z$),
and $\si{s}{m}\in \GG{m}$ is the automorphism sending any $m$th root of unity
to its $s$th power. Hence, an easy observation gives
\begin{equation}\label{Periodicity}
a\equiv b\pmod m\qquad\implies\qquad\theta_m(a)=\theta_m(b).
\end{equation}
Moreover, if $m\mid a$, $\theta_m(a)=0$, whereas if $m\nmid a$ we get the following relation
\begin{equation}\label{Computation_of_theta}
\theta_m(a)+\theta_m(-a)=N_m,
\end{equation} 
where $N_m=\sum_{\tau\in\GG{m}}\tau$ is the absolute norm element.

For any positive integers $m,n$ such that $n\mid m$ we have the usual restriction and corestriction maps between the group rings $\Q[\GG{m}]$ and $\Q[\GG{n}]$
\begin{alignat*}{2}
\res_{\QQ{m}/\QQ{n}}&:\Q[\GG{m}]&&\to\Q[\GG{n}],\\
\cor_{\QQ{m}/\QQ{n}}&:\Q[\GG{n}]&&\to\Q[\GG{m}].
\end{alignat*}
The restriction map is the ring homomorphism sending each automorphism~$\sigma\in\GG{m}$ to its
restriction $\sigma|_{\QQ{n}}$; the corestriction map is the linear map
determined for any $\tau \in\GG{n}$ by
\begin{equation*}
\cor_{\QQ{m}/\QQ{n}}(\tau)=\sum_{\substack{\sigma\in \GG{m}\\\sigma|_{\QQ{n}}=\tau}}\sigma.
\end{equation*}
Let $\Stp m$ be the subgroup of the additive group of $\Q[\GG{m}]$ generated by%
\begin{equation*}
\theta^{(m)}_n(a)=\cor_{\QQ{m}/\QQ{(m,n)}}\Bigl(\res_{\QQ{n}/\QQ{(m,n)}}\bigl(\theta_n(a)\bigr)\Bigr), \text{ for all $a,n\in\Z$, $n>0$.}
\end{equation*}
In fact, $\Stp m$ is Sinnott's group $S'$ from \cite[\crefext{page}{189}]{S80}, for the abelian field~$k$ being the cyclotomic field $\QQ{m}$.%
\footnote{For clarity, let us mention that $\Stp m$ is slightly different from Sinnott's group $S'$ from \cite{S78}.
For example, as $\QQ2=\Q$, the group $\Stp m$ contains $\theta^{(m)}_2(1)=\cor_{\QQ{m}/\Q}\frac12=\frac12N_m$ for each $m$, but Sinnott's group $S'$ from \cite{S78} contains $\frac12N_m$ if and only if $m$ is even.}
The intersection~$\St m=\Stp m\cap\Z[\GG{m}]$ is called the \emph{Stickelberger ideal} of $\QQ{m}$.

\begin{lem}\label{The_index}
For any integer $m>1$, $m\not\equiv2\pmod4$, the index $w=[\Stp m:\St m]$ is equal to the number of roots of unity in the $m$th cyclotomic field $\QQ{m}$, i.e.,
$w=2m$ if $m$ is odd, and $w=m$ if $m$ is even.
\end{lem}
\begin{proof}
This is a part of \cite[\crefext{prop}{2.1}]{S80}.
\end{proof}

\begin{lem}\label{Old_remark}
For any positive integer $m$, the group $\Stp m$ is the subgroup of $\Q[\GG{m}]$ generated by
\begin{equation*}
\bigl\{\theta_m(a);0<a<m\bigr\}\cup\bigl\{\tfrac12N_m\bigr\}.
\end{equation*}
\end{lem}
\begin{proof}
On one hand, $\theta_m(a)=\theta^{(m)}_m(a)\in\Stp m$. On the other hand, let us consider any positive~$n\neq m$ and let $d=(m,n)$. For any $a\in\Z$, using \cite[\crefext{lem}{12}]{K96},
\begin{equation*}
\res_{\QQ{n}/\QQ{d}}\bigl(\theta_n(a)\bigr)\in\Bigl\langle\bigl\{\theta_d(b); 0<b<d\bigr\}\cup\bigl\{\tfrac12N_d\bigr\}\Bigr\rangle.
\end{equation*}
It is easy to see that $\cor_{\QQ{m}/\QQ{d}}\bigl(\tfrac12N_d\bigr)=\tfrac12N_m$.
Considering $\theta_d(b)$,~$0<b<d$,
\begin{equation}\label{Computation}
    \begin{split}
        \cor_{\QQ{m}/\QQ{d}}\bigl(\theta_d(b)\bigr)
        &=\cor_{\QQ{m}/\QQ{d}}\biggl(\ \sum_{\substack{0<s\le d\\(s,d)=1}}\bigl\langle-\tfrac{bs}d\bigr\rangle\si{s}{d}^{-1}\ \biggr)\\&
        =\sum_{\substack{0<s\le m\\(s,m)=1}}\bigl\langle-\tfrac{bs}d\bigr\rangle\si{s}{m}^{-1}
        % \\&
        =\theta_m\bigl(\tfrac{bm}{d}\bigr).
    \end{split}
\end{equation}
As $\cor_{\QQ{m}/\QQ{d}}$ is a group homomorphism, this shows that
\begin{equation*}
\theta^{(m)}_n(a)=\cor_{\QQ{m}/\QQ{d}}\Bigl( \res_{\QQ{n}/\QQ{d}}\bigl(\theta_n(a)\bigr)\Bigr)
\in\Bigl\langle\bigl\{\theta_m(a);0<a<m\bigr\}\cup\bigl\{\tfrac12N_m\bigr\}\Bigr\rangle.
\end{equation*}
The lemma follows.
\end{proof}

We now introduce auxiliary elements that allow to write relations that are useful for the proof of \cref{Old_basisa}.
For any $a\in\Z$, we set
\begin{equation}\label{Definition}
\omega_m(a)=\begin{cases}\theta_m(a)-\frac12N_m,&\text{if }m\nmid a,\\0,&\text{if }m\mid a.\end{cases}
\end{equation}
Adapting \cref{Periodicity,Computation_of_theta,Computation}, we deduce respectively, for $d\mid m$ and $0<b<d$,
\begin{gather}\label{Action_of_-1}
    \omega_m(a+m)=\omega_m(a)
    \qquad\text{and}\qquad
    \omega_m(-a)=-\omega_m(a),\\
    \label{Computation_omega}
    \cor_{\QQ{m}/\QQ{d}}\bigl(\omega_d(b)\bigr) = \cor_{\QQ{m}/\QQ{d}}\bigl(\theta_d(b)-\tfrac12N_d\bigr) =\omega_m\bigl(\tfrac{bm}{d}\bigr).
\end{gather}
The last equality uses that $\cor_{\QQ{m}/\QQ{d}}$ is a linear map and $\cor_{\QQ{m}/\QQ{d}}\bigl(N_d\bigr) = N_m$.
Moreover, by \cref{Old_remark}, $\Stp m$ is the subgroup of $\Q[\GG{m}]$ generated by
\begin{equation}\label{GenSet_omega}
\bigl\{\omega_m(a);0<a<m\bigr\}\cup\bigl\{\tfrac12N_m\bigr\}.
\end{equation}

\begin{lem}\label{Relation}
Let $d,r$ be positive integers and $m=rd$. Then for any $k\in\Z$ we have
\begin{equation*}
\sum_{\substack{a=0,\dots,m-1\\a\equiv k\pmod r}}\omega_m(a)=\sum_{i=0}^{d-1}\omega_m(k+ir)=\omega_m(kd).
\end{equation*}
\end{lem}
\begin{proof}
The lemma follows from the following well-known identity
\begin{equation*}
\sum_{i=0}^{d-1}\left\langle-\frac{s(k+ir)}m\right\rangle=
\left\langle-\frac{skd}m\right\rangle+\frac{d-1}2,
\end{equation*}
valid for any $s\in\Z$ relatively prime to $m$.
\end{proof}

From now on, we shall suppose $m > 1$ is a positive integer, $m\not\equiv 2\pmod 4$.
Let $m=q_1q_2\dots q_t$, where $q_1,q_2,\dots, q_t$ are pairwise coprime prime powers that all satisfy $q_i>2$, and let $p_i$ be the prime dividing $q_i$ for each $i\in\range{1}{t}$.

\begin{xrem}
    Note that we implicitly fix an ordering on the factors~$q_i$ of~$m$. All our results hold true for any ordering as long as it stays consistent through all subsets of the~$q_i$'s. However, if this ambiguity were a problem in an application, we could simply fix an ordering by the assumption $p_1<\dotsb<p_t$.
\end{xrem}

Let $X_m$ be the set of all positive integers $a<m$ that are either divisible by~$q_i$ or relatively prime to $q_i$ for each $i\in\range{1}{t}$, \ie
\begin{equation*}
X_m=\Bigl\{a\in\Z; 0<a<m, \bigl(a,\tfrac m{(a,m)}\bigr)=1\Bigr\}.
\end{equation*}

Let $\ell_i\in\Z$ satisfy $p_i\ell_i\equiv1\pmod{\tfrac{m}{q_i}}$, and~$\ell_i\equiv1\pmod{q_i}$. \Cref{Relation} implies the following result:

\begin{lem} \label{Modified_relation}
For the chosen $m$, for any $i\in\range{1}{t}$ and any $a\in X_m$, we have
\begin{equation*}
\sum_{\substack{k\equiv1\pmod{m/q_i}\\0<k\le m,\ p_i\nmid k}}\omega_m(ka)=
\begin{cases}
\varphi(q_i)\cdot\omega_m(a),&\text{if $q_i\mid a$},\\
\omega_m(aq_i)-\omega_m(aq_i\ell_i),&\text{if $q_i\nmid a$},
\end{cases}
\end{equation*}
where $\varphi$ is Euler's totient function.
\end{lem}

\section{On bases of $\Stp m$}

Recall that $m>1$ is a positive integer such that $m=q_1q_2\dotsc q_t\not\equiv 2\pmod4$, where $\nvars{q}{t}$ are pairwise coprime prime powers greater than $2$.

\subsection{A first basis of $\Stp m$}
We first give a basis of $\Stp m$ constructed in the spirit of \cite[\crefext{theorem}{4.2}]{K92}.
We shall define a useful subset~$M_m$ of the set~$X_m$ defined in the previous section.
Let~$M_m\subseteq X_m$ be the set of all~$a\in X_m$ satisfying
\begin{itemize}
\item for all $i\in\range{1}{t}$, if $q_i\nmid a$ then $a\not\equiv-(a,m)\pmod{q_i}$, 
\item if $a\nmid m$ and 
$k=\max\bigl\{i\in\range{1}{t}; a\not\equiv(a,m)\pmod{q_i}\bigr\}$ then $\bigl\langle\tfrac{a}{(a,m)q_k}\bigr\rangle<\tfrac12$,
\item if $a\mid m$ then the set $\bigl\{ i\in\range{1}{t}; q_i\nmid a\bigr\}$ has an odd number of elements.
\end{itemize}
Actually, $M_m $ is exactly the set $M_{-}$ defined in \cite[\crefext{page}{293}]{K92}. This set has the following stability property:

\begin{lem} \label{Relations_between_sets_M}
Let $r\mid m$, $0<r<m$, such that $\bigl(r,\tfrac mr\bigr)=1$. Let the set $M_{\frac mr}$ be defined using the 
ordering of prime power divisors of $\tfrac mr$ induced by the chosen ordering of prime power divisors of $m$. Then
\begin{equation*}
\bigl\{a\in M_m;\ r\mid a\bigr\}=\bigl\{rb;\ b\in M_{\frac mr}\bigr\}=r\cdot M_{\frac mr}.
\end{equation*}
\end{lem}
\begin{proof}
    For any integer $b$, $0< b < \tfrac mr$, we have $b\in X_{\frac mr}$ if and only if for each $i\in\range{1}{t}$ such that $q_i \mid \tfrac mr$, either $(q_i,b)=1$ or $q_i \mid b$. This is the case if and only if for each $i\in\range{1}{t}$, either $(q_i, rb)=1$ or $q_i \mid rb$, thus if and only if $rb \in X_m$.

If $q_i\mid\frac mr$ for some $i\in\range{1}{t}$ then $(q_i,r)=1$, and so 
$q_i\nmid rb$ if and only if $q_i\nmid b$, moreover $b\not\equiv-(b,\frac mr)\pmod{q_i}$ if and only if $br\not\equiv-(br,m)\pmod{q_i}$.

If $b\nmid\frac mr$ then for any $i\in\range{1}{t}$ such that $q_i\mid\frac mr$ we have $b\not\equiv(b,\frac mr)\pmod{q_i}$ if and only if $br\not\equiv(br,m)\pmod{q_i}$. Therefore we get the same $k$ for $b\in X_{\frac mr}$ and for $br\in X_m$. Moreover $\frac{b}{(b,\frac mr)q_k}=\frac{br}{(br,m)q_k}$.
 
If $b\mid\frac mr$ then $\bigl\{ i\in\range{1}{t}; q_i\mid\frac mr, q_i\nmid b\}=\bigl\{ i\in\range{1}{t}; q_i\nmid rb\}$.
\end{proof}

\begin{thm}\label{Old_basisa}
For any integer $m>1$, $m\not\equiv2\pmod4$, the set 
\begin{equation}\label{Generators}
\bigl\{\omega_m(a);a\in M_m\bigr\}\cup\bigl\{\tfrac12N_m\bigr\}
\end{equation}  
is a $\Z$-basis of $\Stp m$.
\end{thm}

\begin{proof}
This can be proved similarly to the part of \cite[\crefext{thm}{4.2}]{K92} about the Stickelberger ideal, using \cref{Old_remark,Relation,Modified_relation} instead of \cite[\crefext{lem}{3.1},\,3.2, and\,3.4]{K92}. 
Indeed, the proof of \cite[\crefext{thm}{4.2}]{K92} about the Stickelberger ideal and its 
preparatory statements \cite[\crefext{lem}{3.3} and\,4.1 (for $\Psi$)]{K92} 
need the validity of
only the following facts (using notations $\omega(a)$ and $\omega^*$ from \cite{K92}):
\begin{itemize}
    \item the Stickelberger ideal is generated by
$\bigl\{\omega(a);\ 0<a<m\bigr\}\cup\bigl\{\omega^*\bigr\}$ as a group (\cite[\crefext{lem}{3.1}]{K92});
\item these generators satisfy the relations of \cref{Relation}, where we write $\omega(a)$ instead of $\omega_m(a)$
(\cite[\crefext{lem}{3.2}]{K92});
\item these generators satisfy the relations of \cref{Modified_relation}, where we write $\omega(a)$ instead of $\omega_m(a)$
(\cite[\crefext{lem}{3.4}]{K92}).
\end{itemize}
Therefore, this proof can be used \textit{mutatis mutandis} to get a basis for any group generated by generators satisfying these relations. Hence, plugging $\omega^* = \tfrac12 N_m$ and $\omega(a)=\omega_m(a)$, we deduce the theorem from \cref{Old_remark,Relation,Modified_relation}.
\end{proof}

The above basis inherits the stability property given in \cref{Relations_between_sets_M}.

\begin{prop}\label{Old_basisb}
    For any given~$b\in\Z$,~$0<b<m$, let~$r_b$ be the maximal divisor of~$(b,m)$ satisfying~$\bigl(r_b,\tfrac m{r_b}\bigr) = 1$, \ie $r_b$ is the product of all $q_i, i\in\range{1}{t}$ which divide $b$,
    and write $\omega_m(b)\in\Stp m$ as a unique $\Z$-linear combination of basis elements~\eqref{Generators}. %
    Then for each $a\in M_m$ such that $r_b\nmid a$, the coefficient of $\omega_m(a)$ in this $\Z$-linear combination is equal to zero.
\end{prop}

\begin{proof}
    For brevity's sake, let $r=r_b$. By \cref{Computation_omega}, $\omega_m(b)= \cor_{\QQ{m}/\QQ{\frac mr}}\bigl(\omega_{\frac mr}(\tfrac br)\bigr)$. Using \cref{Old_basisa} for $\frac mr$ implies $\omega_{\frac mr}(\frac br)\in\Stp{\frac mr}$ is a unique $\Z$-linear combination of
\begin{equation*}
\bigl\{\omega_{\frac mr}(a);a\in M_{\frac mr}\bigr\}\cup\bigl\{\tfrac12N_{\frac mr}\bigr\}.
\end{equation*}
Since by \cref{Computation_omega}~$\cor_{\QQ{m}/\QQ{\frac mr}}\bigl(\omega_{\frac mr}(a)\bigr)=\omega_{m}(ra)$ and $\cor_{\QQ{m}/\QQ{\frac mr}}\bigl(N_{\frac mr}\bigr)=N_{m}$, and since~$\cor_{\QQ{m}/\QQ{\frac mr}}$ is a linear map, the proposition follows from \cref{Relations_between_sets_M}.
\end{proof}

In particular, for any positive $r\mid m$, $1<r<m$, such that $\bigl(r,\tfrac mr\bigr)=1$, the corestriction subgroup $\cor_{\QQ{m}/\QQ{m/r}}\bigl(\Stp{m/r}\bigr)$ of $\Stp m$ has the following $\Z$-basis
\begin{equation*}
    \bigl\{\omega_m(a); a\in r\cdot M_{\frac mr}\bigr\} \cup \bigl\{\tfrac12 N_m\bigr\}.
\end{equation*}

\subsection{An alternative basis of $\Stp m$: the prime-power case}

In this section we shall suppose that $m$ is a prime power $q=p^e$, where $p$ is a prime and $e$ is a positive integer. Let us mention explicitly that the case $p=2$ is allowed whenever~$e\ge2$ to ensure $q\not\equiv 2\pmod{4}$. We set
\begin{equation}
    M'_q = M'_{p^e}=\bigl\{1,\dots,\tfrac{\varphi(p^e)}2\bigr\}.
\end{equation}

\begin{thm}\label{Prime-power}
For any prime power $q=p^e>2$, the set
\begin{equation}\label{Generators1}
\bigl\{\omega_q(a);a\in M'_q\bigr\}\cup\bigl\{\tfrac12N_q\bigr\}
\end{equation}
is a $\Z$-basis of $\Stp q$.
\end{thm}

\begin{proof}
We shall prove the theorem by induction with respect to $e$. If $q$ is an odd prime or $q=4$, we have $M'_q = M_q$ so this is just a special case of \cref{Old_basisa}.

Let us suppose that the theorem has been proved for $p^e>2$ and let us prove it for $q=p^{e+1}$. 
Let $H$ be the subgroup of $\Stp q$ generated by the set\,\eqref{Generators1}.
We shall show that $H$ contains all $\omega_q(a)$, $0<a<q$, so that it generates $\Stp q$ by \cref{GenSet_omega}.

Since $\omega_q(q-a)=-\omega_q(a)$ by \cref{Action_of_-1},~$H$ contains also~$\omega_q(a)$ for each $a\in\Z$ satisfying $q-\frac{\varphi(q)}2\le a< q$.
Suppose $a=bp$, using \cref{Computation_omega} we get
\begin{equation*}
\omega_q(a) = \cor_{\QQ{q}/\QQ{q/p}}\bigl(\omega_{q/p}(b)\bigr).
\end{equation*}
Since $\cor_{\QQ{q}/\QQ{q/p}}$ is an injective linear map, the induction hypothesis implies that~$\omega_q(a)$ is a linear combination of $\tfrac 12 N_q$ and of $\cor_{\QQ{q}/\QQ{q/p}}\bigl( \omega_{q/p}(t)\bigr) = \omega_q(tp)$ for~$t\in M'_{q/p}$, which implies $tp \in M'_q$.
Thus,~$H$ contains $\omega_q(a)$ whenever $p\mid a$. As for the remaining cases, let $a\in\Z$ such that $\frac{\varphi(q)}2<a<q-\frac{\varphi(q)}2$ and $p\nmid a$.
\Cref{Relation} states that
\begin{equation*}
    \sum_{\substack{t=0,\dots,q-1\\t\equiv a\pmod {q/p}}}\omega_q(t)=\omega_q(ap)\in H.
\end{equation*}
Since $\bigl(q-\tfrac{\varphi(q)}2\bigr)-\tfrac{\varphi(q)}2=\tfrac qp$, there is only one $t$ in the sum on the left hand side satisfying $\tfrac{\varphi(q)}2<t<q-\tfrac{\varphi(q)}2$, namely $t=a$. All other summands are known to belong to~$H$, and since we just proved that $\omega_q(ap)\in H$, we deduce $\omega_q(a)\in H$. 

We have shown that $H$ generates $\Stp q$. Since $|M'_q|=|M_q|$, the theorem follows.
\end{proof}

\subsection{An alternative basis of $\Stp m$: the general case}

Now, we return to the general case where $m = q_1q_2\dotsc q_t \not\equiv 2\pmod4$.
Let us fix~$i\in\range{1}{t}$. \Cref{Relations_between_sets_M} gives that
\begin{equation*}
\bigl\{a\in M_m;\ \tfrac m{q_i}\mid a\bigr\}= \tfrac{m}{q_i}\cdot M_{q_i} = \bigl\{\tfrac {mb}{q_i};\ p_i\nmid b, 0<b<\tfrac{q_i}2 \bigr\}.
\end{equation*}
Since $\cor_{\QQ{m}/\QQ{q_i}}$ is an injective linear map, \cref{Old_basisb} and respectively \cref{Prime-power} combined with \cref{Computation_omega} imply that the sets
\begin{equation*}
    \bigl\{\omega_m(a);\ a\in \tfrac m{q_i}\cdot M_{q_i}\bigr\}\cup\bigl\{\tfrac12N_m\bigr\}
\end{equation*}
and
\begin{equation*}
    \bigl\{\omega_m(b);\ b\in \tfrac{m}{q_i} \cdot M'_{q_i}\bigr\}\cup\bigl\{\tfrac12N_m\bigr\}
\end{equation*}
are $\Z$-bases of the same subgroup $\cor_{\QQ{m}/\QQ{q_i}}\bigl(\Stp{q_i}\bigr)$ of $\Stp m$, so that there is an integral transition matrix between these bases of determinant $\pm1$. We stress that the sets $\tfrac{m}{q_i}\cdot M_{q_i}$ (resp. $\tfrac{m}{q_i}\cdot M'_{q_i}$) for $i\in\range{1}{t}$ are pairwise disjoint. Hence, it is natural to define
\begin{equation}
    \begin{split}
M_m'&=\Bigl(M_m \setminus \bigcup_{i=1}^{t} \tfrac{m}{q_i}\cdot M_{q_i}\Bigr) \cup \Bigl(\bigcup_{i=1}^{t} \tfrac{m}{q_i}\cdot M'_{q_i}\Bigr)\\
&=\Bigl\{a\in M_m;\ \forall i\in\range{1}{t}, \tfrac m{q_i}\nmid a\Bigr\}\cup\Bigl(\bigcup_{i=1}^t
\Bigl\{\tfrac {mb}{q_i};\ 1\leq b\leq {\tfrac{\varphi(q_i)}2}\Bigr\}\Bigr).
\end{split}
\end{equation}
which agrees with the previous definition of $M'_{p^e}$. 
Easily adapting the proof of \cref{Relations_between_sets_M} gives that for any $r\mid m$, $0<r<m$, such that $\bigl(r,\tfrac mr\bigr) = 1$, we have
\begin{equation*}
\bigl\{a\in M'_m;\ r\mid a\bigr\}=\bigl\{rb;\ b\in M'_{\frac mr}\bigr\}=r\cdot M'_{\frac mr}.
\end{equation*}    
Thus, we have proved that \cref{Old_basisa,Old_basisb} implies the following:
\begin{thm}\label{General_casea_omega}
For any integer $m>1$, $m\not\equiv2\pmod4$, the set
\begin{equation}\label{Generators2_omega}
\bigl\{\omega_m(a);\ a\in M_m'\bigr\}\cup\bigl\{\tfrac12N_m\bigr\}
\end{equation}
is a $\Z$-basis of $\Stp m$.
\end{thm}

\begin{prop}
    For any given~$b\in\Z$,~$0<b<m$, let~$r_b$ be the maximal divisor of~$(b,m)$ satisfying~$\bigl(r_b,\tfrac m{r_b}\bigr) = 1$, \ie $r_b$ is the product of all $q_i, i\in\range{1}{t}$ which divide $b$,
    and write $\omega_m(b)\in\Stp m$ as a unique $\Z$-linear combination of basis elements~\eqref{Generators2_omega}. %
Then, for each $a\in M'_m$ such that $r_b\nmid a$, the coefficient of $\omega_m(a)$ in this $\Z$-linear combination is equal to zero.
\end{prop}

Finally, keeping in mind that $\omega_m(a)=\theta_m(a)-\tfrac12N_m$ if $m\nmid a$, we stress that all results of this whole section are equally valid when replacing $\omega_m(\cdot)$ by $\theta_m(\cdot)$, for example:
\begin{corol}\label{General_casea}
For any integer $m>1$, $m\not\equiv2\pmod4$, the set
\begin{equation}\label{Generators2}
\bigl\{\theta_m(a);a\in M_m'\bigr\}\cup\bigl\{\tfrac12N_m\bigr\}
\end{equation}
is a $\Z$-basis of $\Stp m$.
\end{corol}

\begin{corol}\label{General_caseb}
    For any given~$b\in\Z$,~$0<b<m$, let~$r_b$ be the maximal divisor of~$(b,m)$ satisfying~$\bigl(r_b,\tfrac m{r_b}\bigr) = 1$, \ie $r_b$ is the product of all $q_i, i\in\range{1}{t}$ which divide $b$,
    and write $\theta_m(b)\in\Stp m$ as a unique $\Z$-linear combination of basis elements~\eqref{Generators2}. %
Then, for each $a\in M'_m$ such that $r_b\nmid a$, the coefficient of $\theta_m(a)$ in this $\Z$-linear combination is equal to zero.
\end{corol}

\section{Short basis of the Stickelberger ideal}
Elements of $\Z[\GG{m}]$ are called \emph{short} if they are of the form
\begin{equation*}
\sum_{\sigma\in\GG{m}} a_{\sigma} \sigma \in \Z\bigl[\GG{m}\bigr],\qquad\text{where $a_{\sigma}\in\{0,1\}$ for all~$\sigma$.}
\end{equation*}
We first exhibit a large family of short elements of $\St m$. Choosing carefully elements from this family yields a basis\,\eqref{Generators3new} of $\Stp m$ with almost only short elements and also our short basis\,\eqref{Generators6} of the Stickelberger ideal $\St m = \Stp m \cap \Z\bigl[\GG{m}\bigr]$.

\subsection{A family of short elements of $\St m$}
In this section, we construct numerous short elements of $\St m \subset \Stp m$ which we shall use later on.

\begin{prop}\label{Short} 
Let $a,b,c\in\Z$ satisfy $m\nmid a$, $m\nmid b$, $m\nmid c$, $m\mid a+b+c$. Then
\begin{equation*}
\alpha=\theta_m(a)+\theta_m(b)+\theta_m(c)-N_m
\end{equation*}
is a short element of $\St m$. Moreover $(1+\si{-1}{m})\alpha=N_m$, so exactly one half of the coefficients of $\alpha$ are zeros.
\end{prop}

\begin{proof}
Using $\theta_m(c)+\theta_m(-c) = N_m$ when $m \nmid c$ (see \cref{Computation_of_theta}), we obtain
\begin{equation*}
\alpha=
\theta_m(a)+\theta_m(b)-\theta_m(-c)
=\sum_{\substack{0<s\le m\\(s,m)=1}}\Bigl(
\bigl\langle-\tfrac{as}m\bigr\rangle+\bigl\langle-\tfrac{bs}m\bigr\rangle-\bigl\langle\tfrac{cs}m\bigr\rangle
\Bigr)\si{s}{m}^{-1}.
\end{equation*}
Since $0\le\langle x\rangle<1$, every coefficient in the above sum is trivially bounded by
\begin{equation*}
-1<\bigl\langle-\tfrac{as}m\bigr\rangle+\bigl\langle-\tfrac{bs}m\bigr\rangle-\bigl\langle\tfrac{cs}m\bigr\rangle
<2.
\end{equation*}
Moreover, let $[x]=x-\langle x\rangle\in\Z$ be the integral part of $x$ for any $x\in\Q$. Then,
\begin{equation*}
\bigl\langle-\tfrac{as}m\bigr\rangle+\bigl\langle-\tfrac{bs}m\bigr\rangle-\bigl\langle\tfrac{cs}m\bigr\rangle=
-\tfrac{(a+b+c)s}m
-\bigl[-\tfrac{as}m\bigr]-\bigl[-\tfrac{bs}m\bigr]+\bigl[\tfrac{cs}m\bigr]\in\Z,
\end{equation*}
which proves that $\alpha$ is short. The last equality of the proposition follows again from \cref{Computation_of_theta} and an easy observation that $\si{-1}{m}\theta_m(a)=\theta_m(-a)$.
\end{proof}

\subsection{Bases of $\Stp m$ with many short elements}\label{Definition_of_alpha}
We first describe the map $\alpha_m$, which associates to any $b\in\Z$, $0<b<m$, one short element from the family of \cref{Short}.
For any given $b\in\Z$, let $r_b$ be the maximal divisor $r$ of $(b,m)$ satisfying $(r,\frac m{r})=1$. In other words,
\begin{equation*}
r_b=\prod_{i\in J_b}q_i, \qquad\text{where $J_b=\bigl\{i\in\range{1}{t}; q_i\mid b\bigr\}$.}
\end{equation*}
Let $J_b'=\range{1}{t}\setminus J_b=\bigl\{i\in\range{1}{t}; q_i\nmid b\bigr\}$, and let us suppose that $0<b<m$ so that $J_b'\ne\emptyset$. We define $\alpha_m(b)$ as follows:
\begin{itemize}
\item If $|J_b'|>1$, let $u=q_{\min J_b'}$, and $v=\frac{m}{ur_b}$. 
Since $(u,v)=1$, the equation
\begin{equation*}
ux+vy=-1
\end{equation*}
has a solution $x,y\in\Z$, where $x$ is well-defined modulo $v$ and $y$ modulo $u$,
so $bux$ and $bvy$ are well-defined modulo $m$.
Let
\begin{equation}\label{First_transition}
\alpha_m(b)=\theta_m(b)+\theta_m(bux)+\theta_m(bvy)-N_m.
\end{equation}
\item If $J_b'=\{j\}$ then $b=\frac{mc}{q_j}$ for a unique $c\in\Z$, $0<c<q_j$. If $c>1$ we define
\begin{equation}\label{Second_transition}
\alpha_m(b)=\theta_m(-b)+\theta_m\bigl(b-\tfrac m{q_j}\bigr)+\theta_m\bigl(\tfrac m{q_j}\bigr)-N_m,
\end{equation}
whereas if $c=1$, so that $b=\tfrac m{q_j}$, we put
\begin{equation}\label{Third_transition}
\alpha_m(b)=2\theta_m\bigl(\tfrac{m\cdot\varphi(q_j)}{2q_j}\bigr)+\theta_m\bigl(\tfrac m{p_j}\bigr)-N_m.
\end{equation}
\end{itemize}

Intuitively, $\alpha_m(\cdot)$ is constructed by means of layers on $|J'_b|$, similarly to what happens for $M_m$ as shown by \cref{Relations_between_sets_M}. For $|J'_b|=1$, we follow the prime power case of \cref{Prime-power}, which is very similar to \cite[\crefext{thm}{9.3(i)}]{Sch10} when $m=p$. For~$|J'_b|>1$ we use Bezout's equality to write $-b$ as the sum of two summands~$bux$ and~$buy$ in such a way that both $|J'_{bux}|$ and $|J'_{bvy}|$ are strictly smaller than $|J'_b|$, so that both~$\theta_m(bux)$ and~$\theta_m(bvy)$ are generated by basis elements that were already chosen in the previous layers. Any way of achieving this property works.
In particular, note that in the case $|J_b'|>1$ we could use any other decomposition of $\frac{m}{r_b}$ into the product of relatively prime integers $u>1$, $v>1$.

\begin{lem}\label{Alpha_is_short}
For any integer $m>1$, $m\not\equiv2\pmod4$, the element $\alpha_m(b)$ is short and satisfies $(1+\si{-1}{m})\alpha_m(b)=N_m$ for each positive integer $b<m$.
\end{lem} 
\begin{proof}
In the former case $|J_b'|>1$, we have~$b+bux+bvy=0$. Since  $u\nmid b$, we have~$u\nmid bvy$; similarly $v\nmid b$ implies $v\nmid bux$. Hence~$\alpha_m(b)$ is short by \cref{Short}.
In the latter case~$J_b'=\{j\}$ for some $j\in\range{1}{t}$, we have that~$b$ writes as~$\frac{mc}{q_j}$ with~$c\in\Z$ and~$0<c<q_j$, then $\alpha_m(b)$ is short by \cref{Short} again, because~$-b+(b-\tfrac m{q_j})+\tfrac m{q_j}=0$ and~$2\cdot \tfrac{m\cdot\varphi(q_j)}{2q_j}+\tfrac m{p_j}=m$.
\end{proof}

\begin{thm}\label{Almost_short}
For any integer $m>1$, $m\not\equiv2\pmod4$, the sets
\begin{align}\label{Generators3}
\Bigl\{\alpha_m(b);b\in M_m',|J_b'|>1\Bigr\}\cup\Bigl\{\theta_m(b);b\in M_m',|J_b'|=1\Bigr\}\cup\Bigl\{\tfrac12N_m\Bigr\},\\
  \label{Generators3new}
\quad\Bigl\{\alpha_m(b);b\in M_m'\setminus\bigl\{\tfrac m{q_1},\dotsc,\tfrac m{q_t}\bigr\}\Bigr\}\cup\Bigl\{
\theta_m\bigl(\tfrac m{q_1}\bigr),\dotsc,\theta_m\bigl(\tfrac m{q_t}\bigr),\tfrac12N_m\Bigr\}
\end{align}
are $\Z$-bases of $\Stp m$.
\end{thm}

\begin{proof}
By definition of $\alpha_m(b)$ in \cref{First_transition,Second_transition}, we know that all elements of these sets belong to $\Stp m$. We shall show that the transition matrices from the set\,\eqref{Generators2} to the set\,\eqref{Generators3} and from the set\,\eqref{Generators3} to the set\,\eqref{Generators3new} are, after a suitable reordering of elements of $M'_m$, triangular with~$\pm1$ on the diagonal, which will prove the theorem.

At first, we deal with the transition matrix from the set\,\eqref{Generators2} to the set\,\eqref{Generators3} and 
we shall use induction with respect to $|J_b'|$.
If $|J_b'|=1$ then $\theta_m(b)$  belongs to both sets\,\eqref{Generators2} and\,\eqref{Generators3}.
So suppose that $|J_b'|>1$. Then the transition from~$\theta_m(b)$ to $\alpha_m(b)$ given in \cref{First_transition} uses
$\theta_m(bux)$ and $\theta_m(bvy)$ and the coefficient of $\theta_m(b)$ is $1$.
By \cref{General_caseb},~$\theta_m(bux)$ is a $\Z$-linear combination of $\theta_m(a)$ for $a$ running over $M'_m$ such that $r_{bux} \mid a$. For these $a$'s, we have that
\begin{equation*} 
J'_a \subseteq J'_{r_{bux}} = J'_{bux} \subsetneq J'_b,
\end{equation*}
since $\min J'_b\notin J'_{bux}$ by definition of $u$. Hence, all these $\theta_m(a)$ are covered by induction, and so is $\theta_m(bux)$. The case of $\theta_m(bvy)$ can be treated similarly.

Now, let us consider  the transition matrix from the set\,\eqref{Generators3} to the set\,\eqref{Generators3new}.
Suppose that $J_b'=\{j\}$ and $b=\frac{mc}{q_j}$ for some $c\in\Z$, $1\le c\le\tfrac{\varphi(q_j)}2$.
If~$c=1$ then~$\theta_m(b)$ belongs to both sets\,\eqref{Generators3} and\,\eqref{Generators3new}.
If $c>1$ then the transition from~$\theta_m(b)$ to $\alpha_m(b)$, by \cref{Second_transition,Computation_of_theta}, writes as
\begin{equation*}
\alpha_m(b)=-\theta_m(b)+\theta_m\bigl(b-\tfrac m{q_j}\bigr)+\theta_m\bigl(\tfrac m{q_j}\bigr).
\end{equation*}
Since $J_{b-m/q_j}'=J_{m/{q_j}}'=J_b'$, both $\theta_m\bigl(b-\tfrac m{q_j}\bigr)=\theta_m\bigl(\tfrac m{q_j}(c-1)\bigr)$ and $\theta_m\bigl(\tfrac m{q_j}\bigr)$ were already covered by induction. The coefficient of $\theta_m(b)$ is $-1$.
\end{proof}

\subsection{A basis of $\St m$ with only short elements}

Recall that the Stickelberger ideal of $\QQ{m}$ is the intersection $\St m=\Stp m\cap\Z[\GG{m}]$. Let $\Stpp m$ be the subgroup of $\Stp m$ generated by the set
\begin{equation}\label{Generators4}
\bigl\{\alpha_m(a);\ a\in M_m'\bigr\}\cup\bigl\{\tfrac12N_m\bigr\}.
\end{equation}
We shall prove that $\Stpp m=\St m + \tfrac 12 N_m\cdot\Z$ and that~\cref{Generators4} is its basis. We shall start by computing its finite index in $\Stp m$.
First, we treat the prime power case.
\begin{lem}\label{Special_case}
Let $q=p^e>2$, where $p$ is a prime and $e$ is a positive integer. 
Then the index of $\Stpp q$ in $\Stp q$ is finite and
\begin{equation*}
[\Stp q:\Stpp q]=\begin{cases}\frac q2&\text{if $p=2$},\\q&\text{if $p>2$}.
\end{cases}
\end{equation*}
\end{lem}
\begin{proof}
    To obtain the index $[\Stp q:\Stpp q]$, let us compute the transition matrix from
\begin{equation}\label{Generators2_prime-power}
\bigl\{\theta_q(a);\ a\in\Z, 1\le a\le \tfrac{\varphi(q)}2\bigr\}\cup\bigl\{\tfrac12N_q\bigr\},
\end{equation}
which is a $\Z$-basis of $\Stp q$ by \cref{General_casea}, to the system of generators of $\Stpp q$, \ie
\begin{equation}\label{Generators5}
\bigl\{\alpha_q(a);\ a\in\Z, 1\le a\le \tfrac{\varphi(q)}2\bigr\}\cup\bigl\{\tfrac12N_q\bigr\}.
\end{equation}
This transition matrix is given by \cref{Third_transition,Second_transition}.
More precisely, using also \cref{Computation_of_theta}, we obtain in the studied special case that
\begin{equation*}
    \alpha_q(a) = \begin{cases}\theta_q\bigl(p^{e-1}\bigr)+2\theta_q\bigl(\tfrac{\varphi(q)}2\bigr)-N_q&\text{if $a=1$},\\ \theta_q(1)+\theta_q(a-1)-\theta_q(a)&\text{if $2\leq a \leq\frac{\varphi(q)}{2}$}.
    \end{cases}
\end{equation*}
Since $\frac12N_q$ belongs to both sets\,\eqref{Generators2_prime-power} and\,\eqref{Generators5}, we can ignore this element in the computation of the determinant of the transition matrix.

At first, let us assume that $p>3$. Then $p^{e-1} < p^{e-1}\cdot\tfrac{p-1}2=\tfrac{\varphi(q)}2$.
We shall compute the determinant of the following square matrix of dimension $\tfrac{\varphi(q)}2$
\begin{equation}\label{Matrix}
\begin{pmatrix}
0&0&0&0&\cdots &1 &\cdots &0&0&2\\
2&-1&0&0&\cdots &0 &\cdots &0&0&0\\
1&1&-1&0&\cdots &0 &\cdots &0&0&0\\
1&0&1&-1&\cdots &0 &\cdots &0&0&0\\
\vdots &\vdots &\vdots &\vdots & &\vdots & &\vdots &\vdots  &\vdots\\
1&0&0&0&\cdots &0 &\cdots &1&-1&0\\
1&0&0&0&\cdots &0 &\cdots &0&1&-1
\end{pmatrix},
\end{equation}
where the $1$ in the first row belongs to the $p^{e-1}$th column (which is the first column if $e=1$).
The sum of all rows but the first one, multiplied by $2$, equals
\begin{equation*}
\begin{pmatrix}
\varphi(q)&0&0&0&\cdots &0 &\cdots &0&0&-2
\end{pmatrix}.
\end{equation*}
We add this row to the first row of our matrix.
If $e>1$, we also add to the first row the sum of all rows from the second one to the $p^{e-1}$th one, \ie
\begin{equation*}
\begin{pmatrix}
p^{e-1}&0&0&0&\cdots &-1 &\cdots &0&0&0
\end{pmatrix}.
\end{equation*}
After this computation we get a lower triangular matrix of determinant $\pm q$. As this determinant is nonzero, the set\,\eqref{Generators5} is a $\Z$-basis of $\Stpp q$ and the index~$[\Stp q:\Stpp q]$ equals the absolute value of the determinant. The lemma follows for~$p>3$.

Now, suppose $p=3$. Then $p^{e-1} =\tfrac{\varphi(q)}2$
and the square transition matrix of dimension~$3^{e-1}$ writes as
\begin{equation*}
\begin{pmatrix}
0&0&0&0&\cdots &0&0&3\\
2&-1&0&0&\cdots &0&0&0\\
1&1&-1&0&\cdots &0&0&0\\
1&0&1&-1&\cdots &0&0&0\\
\vdots &\vdots &\vdots &\vdots & &\vdots &\vdots  &\vdots\\
1&0&0&0&\cdots &1&-1&0\\
1&0&0&0&\cdots &0&1&-1
\end{pmatrix}.
\end{equation*}
If $e=1$ then the only entry of our matrix of dimension $1$ is $3$. If $e>1$, the
sum of all rows but the first one, multiplied by $3$, is equal to
\begin{equation*}
\begin{pmatrix}
3^e&0&0&0&\cdots &0 &\cdots &0&0&-3
\end{pmatrix}.
\end{equation*}
Adding this row to the first row, we again get a lower triangular matrix of determinant $\pm q$, which gives the lemma in the case $p=3$.

Finally, we treat the case $p=2$. Then, by \cref{Periodicity,Computation_of_theta}, we have
\begin{equation*}
\theta_q(2^{e-1})=\theta_q(2^{e-1}-q)=\theta_q(-2^{e-1})=N_q-\theta_q(2^{e-1}),
\end{equation*}
so $\theta_q(2^{e-1})=\frac12N_q$. Therefore we have got almost the same matrix as written in \cref{Matrix}, except that in the first row the only non-zero element is the 2 at the very end. By the same approach as above, we obtain that the determinant of this matrix is equal to $\pm\varphi(q)=\pm\frac q2$ and the lemma in the case $p=2$ follows.
\end{proof}

\begin{prop}\label{General_index}
For any integer $m>1$, $m\not\equiv2\pmod4$, the set\,\eqref{Generators4} is a basis of $\Stpp m$, whose finite index in $\Stp m$ is given by
\begin{equation*}
[\Stp m:\Stpp m]=\begin{cases}\frac m2&\text{if $m$ is even},\\m&\text{if $m$ is odd}.
\end{cases}
\end{equation*}
\end{prop}
\begin{proof}
This is similar to the proof of \cref{General_casea_omega}.
The following sets are pairwise disjoints for $i \in \range{1}{t}$
\begin{equation*}
\bigl\{\tfrac {mb}{q_i};\ b\in M'_{q_i}\bigr\}=\bigl\{a;\ a\in M'_m,\ \tfrac {m}{q_i}\mid a\bigr\}.
\end{equation*}
Since $\cor_{\QQ{m}/\QQ{q_i}}$ is an injective linear map, 
the transition matrix from the $\Z$-basis\,\eqref{Generators3}
of $\Stp m$, given by \cref{Almost_short}, to the system of generators\,\eqref{Generators4}
of $\Stpp m$ is a block diagonal matrix, having (besides plenty of trivial blocks of dimension~$1$ containing~$1$) one nontrivial block for each $i\in\range{1}{t}$. For a given $i$, the nontrivial block is equal to the matrix considered in \cref{Special_case} for $q=q_i$. Since the determinant of this transition matrix is equal to the product of determinants of these nontrivial blocks, it is nonzero and the proposition follows.  
\end{proof}

We are now ready to state our main theorem, which in particular implies the afore-mentioned relation $\Stpp m = \St m + \tfrac 12 N_m\cdot\Z$.

\begin{thm}\label{Short_basis}
For any integer $m>1$, $m\not\equiv2\pmod4$, the set
\begin{equation}\label{Generators6}
\bigl\{\alpha_m(a);a\in M_m'\bigr\}\cup\bigl\{N_m\bigr\}
\end{equation}
is a $\Z$-basis of the Stickelberger ideal $\St m$ of $\QQ{m}$ having only short elements. 
\end{thm}
\begin{proof}
Let $\widetilde{\St m}$ denote the subgroup of $\Stp m$ generated by the set\,\eqref{Generators6}.
Each element of \eqref{Generators6} is short by \cref{Alpha_is_short}, in particular it belongs to $\Z[\GG{m}]$, so that
\begin{equation}\label{Inclusion}
\widetilde{\St m}\subseteq\Z[\GG{m}]\cap\Stp m=\St m.
\end{equation}
The indices $[\Stp m:\St m]=w$ and $[\Stp m:\Stpp m]=\frac w2$ are given by \cref{The_index,General_index}, respectively.
In particular, by \cref{General_index}, the set\,\eqref{Generators4} is linearly independent;
comparing with the set\,\eqref{Generators6}, we see that the set\,\eqref{Generators6} is also linearly independent and that $\widetilde{\St m}$ is a subgroup of $\Stpp m$ of index 
$[\Stpp m:\widetilde{\St m}]=2$. Hence,
\begin{equation*}
[\Stp m:\widetilde{\St m}]=[\Stp m:\Stpp m]\cdot[\Stpp m:\widetilde{\St m}]=w=[\Stp m:\St m],
\end{equation*}
and the inclusion\,\eqref{Inclusion} gives $\widetilde{\St m}=\St m$.
The theorem follows.
\end{proof}

\section{An upper bound for the relative class number of a cyclotomic field}
\label{Relative_hK}
Our short basis of the Stickelberger ideal $\St m$, given in \cref{Short_basis}, allows to derive a simple upper bound on the relative class number of \emph{any} cyclotomic field.

\begin{corol}\label{Upper_bound}
Let $m>1$ be an integer satisfying $m\not\equiv2\pmod4$, let $t$ be the number of primes dividing $m$. The relative class number $h_m^-$ of the $m$th cyclotomic field satisfies
\begin{equation*}
h_m^-\le2^{1-a} \cdot \bigl(\tfrac{\varphi(m)}8\bigr)^{\varphi(m)/4},
\end{equation*}
where $\varphi$ is Euler's totient function and
\begin{equation}\label{Definition_of_a}
a=\begin{cases}
0&\text{if $t=1$},\\2^{t-2}-1&\text{if $t\ge2$}.
\end{cases}
\end{equation}
\end{corol}
\begin{proof}
    Recall that, for any integer $s$ relatively prime to $m$, $\si{s}{m}\in \GG{m}$ denotes the automorphism of the $m$th cyclotomic field $\QQ{m}$ sending any $m$th root of unity
to its $s$th power. In particular, $\si{-1}{m}$ is the restriction of the complex conjugation.
Following Sinnott, let $\R_m=\Z[\GG{m}]$ and
\begin{align*}
\R_m^-&=\{\alpha\in \R_m; (1+\si{-1}{m})\alpha=0\},\\
\A_m&=\{\alpha\in \R_m; (1+\si{-1}{m})\alpha\in N_m\Z\}.
\end{align*}
Moreover, for any submodule $M\subseteq \R_m$ we define $M^-=M\cap \R_m^-$. Using \cite[\crefext{lem}{1.2(a)}]{S80}, multiplication by $1+\si{-1}{m}$ gives
\begin{equation*}
[\A_m:\St m]=[(1+\si{-1}{m})\A_m:(1+\si{-1}{m})\St m]\cdot[\A_m^-:\St m^-].
\end{equation*}
It is clear that $(1+\si{-1}{m})\A_m=(1+\si{-1}{m})\St m=N_m\Z$ and that $\A_m^-=\R_m^-$. Therefore, using \cite[Th., \crefext{page}{107}]{S78}, we have
\begin{equation}\label{Sinnott_result}
[\A_m:\St m]=[\R_m^-:\St m^-]=2^{a}\cdot h_m^-,
\end{equation}
where $a$ is defined by \cref{Definition_of_a}.

We use our short basis\,\eqref{Generators6} of $\St m$ given in \cref{Short_basis} to derive a bound on $[\A_m:\St m]$. First, a $\Z$-basis of $\A_m$ is given by
\begin{equation}\label{Basis_of_A}
\bigl\{\beta_m(s); 1\le s<\tfrac m2, (s,m)=1\bigr\}\cup\bigl\{\gamma_m\bigr\},
\end{equation}
where $\beta_m(s)=\si{s}{m}-\si{-s}{m}$ and
\begin{equation*}
\gamma_m=\sum_{\substack{1\le s<\frac m2\\(s,m)=1}}\si{s}{m}.
\end{equation*}
An easy calculation gives
\begin{equation*}
N_m=2\gamma_m-\sum_{\substack{1\le s<\frac m2\\(s,m)=1}}\beta_m(s).
\end{equation*}
For each $b\in M_m'$, let us define integers $a_{b,s}$, where $1\le s<m$, $(s,m)=1$, by
\begin{equation*}
\alpha_m(b)=\sum_{\substack{1\le s<m\\(s,m)=1}}a_{b,s}\si{s}{m}.
\end{equation*}
By \cref{Alpha_is_short},
we have $a_{b,s}+a_{b,m-s}=1$, so that
\begin{equation*}
\alpha_m(b)=\gamma_m + \sum_{\substack{1\le s<\frac m2\\(s,m)=1}}(a_{b,s}-1)\beta_m(s).
\end{equation*}
The index $[\A_m:\St m]$ is given by the absolute value of the determinant of the transition matrix from the basis\,\eqref{Generators6} of $\St m$ to the basis\,\eqref{Basis_of_A} of $\A_m$, \ie
\newlength{\hmin}
\newlength{\wmin}
\newlength{\wmax}
\settoheight{\hmin}{$\vdots$}
\settowidth{\wmin}{$-1$}
\settowidth{\wmax}{$a_{b,s}-1$}
\begin{equation*}
    \bigl[\A_m:\St m\bigr] = \left|\det\!\!
    \begin{tikzpicture}[
        baseline=(AmSm.center),
        every left delimiter/.style={xshift=1.1em},
        every right delimiter/.style={xshift=-.2em},
    ]
        \matrix (AmSm) [
            every node/.style={anchor=center},
            matrix of math nodes,
            left delimiter={(}, right delimiter={)},
            nodes in empty cells,
            nodes={minimum height=\hmin+4pt,minimum width=\wmax,inner sep=2pt},
        ] {
          |[minimum width=\wmin]|2      & -1 & \hdots & -1 \\
          |[minimum width=\wmin]|1      &&&\\
          |[minimum width=\wmin]|\vphantom{\int^0}\smash[t]{\vdots} &&&\\
          |[minimum width=\wmin]|1      &&&\\
        };
        \node[fit=(AmSm-2-2)(AmSm-4-4),xshift=-.3\wmax,yshift=.2\hmin] {$\smash{\biggl( a_{b,s}-1\biggr)}$};
        \node[align=center,anchor=center] (truc) at ([xshift=-.8em]AmSm-4-4.north) {$\substack{b\in M'_m\\1\leq s<\frac m2,\,(s,m)=1}$};
        % Optional lines
        \draw ([yshift=0em]AmSm-1-1.north east) -- ([yshift=0em]AmSm-4-1.south east);
        \draw ([xshift=0.3em]AmSm-1-1.south west) -- ([xshift=0.6em]AmSm-1-4.south east);
    \end{tikzpicture}\!\right|.
\end{equation*}
We subtract one half of the first row from each of the other rows to get
\begin{equation*}
    \bigl[\A_m:\St m\bigr] = \left|\det\!\!
    \begin{tikzpicture}[
        baseline=(AmSmPivot.center),
        every left delimiter/.style={xshift=1.1em},
        every right delimiter/.style={xshift=-.2em}
    ]
        \matrix (AmSmPivot) [
            every node/.style={anchor=center},
            matrix of math nodes,
            left delimiter={(}, right delimiter={)},
            nodes in empty cells,
            nodes={minimum height=\hmin+4pt,minimum width=\wmax,inner sep=2pt},
        ] {
          |[minimum width=\wmin]|2      & -1 & \hdots & -1 \\
          |[minimum width=\wmin]|0      &&&\\
          |[minimum width=\wmin]|\vphantom{\int^0}\smash[t]{\vdots} &&&\\
          |[minimum width=\wmin]|0      &&&\\
        };
        \node[fit=(AmSmPivot-2-2)(AmSmPivot-4-4),xshift=-.3\wmax,yshift=.2\hmin] {$\smash{\biggl( a_{b,s}-\tfrac 12\biggr)}$};
        \node[align=center,anchor=center] at ([xshift=-.8em]AmSmPivot-4-4.north) {$\substack{b\in M'_m\\1\leq s<\frac m2,\,(s,m)=1}$};
        % Optional lines
        \draw ([yshift=0em]AmSm-1-1.north east) -- ([yshift=0em]AmSm-4-1.south east);
        \draw ([xshift=0.3em]AmSm-1-1.south west) -- ([xshift=0.6em]AmSm-1-4.south east);
    \end{tikzpicture}\!\right|
    = 2\cdot\left|\det\!\!
    \begin{tikzpicture}[
        baseline=(AmSmFinal.center),
        every left delimiter/.style={xshift=.9em},
        every right delimiter/.style={xshift=.3em},
    ]
        \matrix (AmSmFinal) [
            every node/.style={anchor=center},
            matrix of math nodes,
            left delimiter={(}, right delimiter={)},
            nodes in empty cells,
            nodes={minimum height=\hmin+4pt,minimum width=.8\wmax,inner sep=2pt},
        ] {
          &&\\
          &&\\
          &&\\
        };
        \node[fit=(AmSmFinal-1-1)(AmSmFinal-3-3),align=left,yshift=.2\hmin] {$\smash{\biggl( a_{b,s}-\tfrac 12\biggr)}$};
        \node[align=center,anchor=center] at ([xshift=-.7em]AmSmFinal-3-3.north) {$\substack{b\in M'_m\\1\leq s<\frac m2,\,(s,m)=1}$};
    \end{tikzpicture}\!\right|.
\end{equation*}
By \cref{Alpha_is_short} we know that $a_{b,s}\in\{0,1\}$, and so $a_{b,s}-\tfrac12\in\bigl\{-\tfrac12,\tfrac12\bigr\}$. So the length of each row of this matrix, as a vector in the Euclidean space of dimension $\tfrac{\varphi(m)}{2}$, is equal to $\tfrac12\sqrt{\tfrac{\varphi(m)}2}$.
Therefore, by Hadamard's inequality,
\begin{equation*}
[\A_m:\St m]\le 2\cdot\biggl(\frac12\sqrt{\frac{\varphi(m)}2}\biggr)^{\varphi(m)/2}.
\end{equation*}
A comparison with \cref{Sinnott_result} gives the corollary.
\end{proof}

\begin{xrem}
For the marginal cases where $4\nmid\frac{\varphi(m)}2$, better bounds exist for these scaled Hadamard matrices (see \cite{BEHC21}) that directly translate into slightly better bounds for $h_m^-$. We do not dive into the details here.
\end{xrem}

\section{Effective short Stickelberger generators}
Let $m >1$ satisfy $m \not \equiv 2\pmod4$. Let $\ell$ be any prime such that $(\ell,m) = 1$ and let~$\LL$ be a fixed (unramified) prime ideal above $\ell$ of inertia degree~$f$ in the~$m$th cyclotomic field~$\QQ{m}$. The aim of this section is to describe an algebraic integer of $\QQ{m}$ generating the principal ideal $\LL^{\alpha_m(b)}$ for each $b\in M'_m$.

Of course, we shall use Gauss sums.
Recall that $\zeta_n = e^{2\pi i/ n}$ for any positive integer $n$.
Let $\F = \Z[\zeta_m]\bigl/ \LL$ be the finite field of cardinality $\norm{\LL}=\ell^f$, and
let~$\chi_{\LL}$ be the~$m$th power Legendre symbol with respect to~$\LL$, i.e., for any $a\in\F^{\times}$, $\chi_\LL(a) \in\bigl\langle \zeta_m\bigr\rangle$ is 
determined by the condition that $\chi_\LL(a)$ belongs to the class $a^{ (\norm{\LL}-1)/m }$.
We extend as usual characters to $\F$ by setting $\chi_{\LL}(0) = 0$.
For any integer $b$, we have the following Gauss sum, where $\TraceOp:\F \rightarrow \F_{\ell}$ is the trace map in the field extension $\F \bigl/ \F_{\ell}$,
\begin{equation*}
    g(b,\LL) = -\sum_{a\in\F}  \chi_\LL(a)^b \zeta_{\ell}^{\TraceOp(a)}\in\Z[\zeta_{m\ell}].
\end{equation*}
For any integers $u\equiv1\pmod m$, $\ell\nmid u$, and
$v\equiv1\pmod\ell$, $(v,m)=1$,
an easy computation gives (see \eg \cite[(3.3)~and\,(3.5)]{S80})
\begin{align}\label{Galois_action1}
\si{u}{m\ell}\bigl(g(b,\LL)\bigr)&=\chi_\LL(u)^{-b}\cdot g(b,\LL),\\
\label{Galois_action2}
\si{v}{m\ell}\bigl(g(b,\LL)\bigr)&=g(vb,\LL).
\end{align}
Hence, $g(b,\LL)^m\in\Z[\zeta_m]$ by \cref{Galois_action1}. Moreover, we have the well-known Stickelberger factorization (see \eg \cite[(3.4)]{S80})
\begin{equation}\label{Stickelberger_factorization}
g(b,\LL)^m \cdot\Z[\zeta_m] = \LL^{m\theta_m(b)}.
\end{equation}

We want to describe an explicit generator of the principal ideal 
$\LL^{\alpha_m(b)}$ for each $b\in M'_m$. Since each $\alpha_m(b)$ is given by the general construction from \cref{Short} (see the proof of \cref{Alpha_is_short}), we shall start more generally.

\begin{prop}\label{Jacobi_sum}
For any $b,c\in\Z$ such that $m\nmid b$, $m\nmid c$, $m\nmid b+c$, let
\begin{equation*}
\alpha=\theta_m(b)+\theta_m(c)-\theta_m(b+c)\in\Z[\GG{m}]
\end{equation*}
be one of the short elements given by \cref{Short}.
Then the Jacobi sum
\begin{equation*}
    J(b,c,\LL)=-\sum_{a\in\F} \chi_\LL(a)^b \chi_\LL(1-a)^c \in \Z[\zeta_m]
\end{equation*}
satisfies
$J(b,c,\LL)\cdot\Z[\zeta_m]=\LL^\alpha$.
\end{prop}
\begin{proof}
By \cite[\crefext{lem}{6.2(d)}]{Was97}, we have
\begin{equation*}
J(b,c,\LL)=\frac{g(b,\LL)g(c,\LL)}{g(b+c,\LL)}.
\end{equation*}
Thus, the result directly follows from \cref{Stickelberger_factorization} and the fact $J(b,c,\LL)\in\Z[\zeta_m]$.
\end{proof}

As an example of application of \cref{Jacobi_sum}, let us consider any $b\in M'_m$ such that $|J'_b|>1$.
Then $\alpha_m(b)$ is given by \cref{First_transition}, so that
\begin{equation*}
\LL^{\alpha_m(b)}=J(bux,bvy,\LL)\cdot \Z[\zeta_m], 
\end{equation*}
where $u=q_{\min J_b'}$, $v=\frac{m}{ur_b}$, and $x,y\in\Z$
satisfy $ux+vy=-1$.

Furthermore, it is clear that $u$, $v$, $x$, $y$ do not depend on $b$ but only on
$J'_b$. Therefore, having another $c\in M'_m$ such that $J'_c=J'_b$,
there is an integer $s$ relatively prime to $m$ satisfying 
$c\equiv sb\pmod m$, so that \cref{Galois_action2} gives
\begin{equation*}
J(cux,cvy,\LL)=J(sbux,sbvy,\LL)=\si{s}{m}\bigl(J(bux,bvy,\LL)\bigr).
\end{equation*}
Hence, computing generators for all $\LL^{\alpha_m(b)}$, $b\in M'_m$, comes down to the computation of exactly one representative Jacobi sum per set $J'_b$, then applying a suitable automorphism to obtain the generator for $\LL^{\alpha_m(c)}$ whenever $J'_c = J'_b$.

\newcommand{\sage}{\textsc{SageMath}}
\newcommand{\cpu}{Intel\textsuperscript{\textregistered} Core{\texttrademark} i7-8650U @3.2GHz}

\section{Practical results}

We implemented in practice the computation of our short Stickelberger bases from \cref{Short_basis} using \sage{} \cite{Sage} on an \cpu{}.

All involved algebraic criteria are very easy to compute, so that obtaining the short bases is actually a matter of seconds for any reasonable conductor.
We verified, for all conductors $m<10000$, $m\not\equiv 2\pmod4$, such that $\varphi(m)\leq 2000$, that the Hermite Normal Form (HNF) of the short basis from \cref{Short_basis} coincides with the HNF of the large basis from \cite[\crefext{thm}{6.2}]{K92}.

We stress that using a naive trial-and-error strategy to extract a short basis from a large set of short vectors, \eg from the set $W$ of \cite[\crefext{section}{4.2}]{CDW21}, may converge only after a huge number of iterations, each involving the computation of a costly HNF. This is especially hazardous when $t$ grows, \eg our brute force experiment for $m=780=2^2\cdot3\cdot5\cdot13$ never finished despite the small dimension.

% More interestingly: Computation of h- using [Am:Sm]
More interestingly, we used the determinant formula for $\bigl[\A_m:\St m\bigr]$ given in~\cref{Relative_hK} to derive the relative class number $h_m^-$ from \cref{Sinnott_result}. We checked, for the same range of conductors as above, that the obtained values coincide with the values given by the analytic class number formula \cite[\crefext{thm}{4.17}]{Was97}
\begin{equation}\label{analytic_hminus}
    h_m^- = Qw\cdot\prod_{\text{$\chi$ odd}} \Bigl(-\tfrac12 B_{1,\chi}\Bigr),
\end{equation}
% where $w=2m$ if $m$ is odd and $w=m$ if $m$ is even, $Q=1$ if $m$ is a prime power and $Q=2$ otherwise, and $B_{1,\chi}$ is defined by $\tfrac1f \sum_{a=1}^f a\cdot\chi(a)$ for any odd primitive character $\chi$ modulo $m$ of conductor~$f$.
where the product is taken over all odd primitive Dirichlet characters $\chi$ of conductor $f_{\chi}\mid m$, $w=2m$ if $m$ is odd and $w=m$ if $m$ is even, $Q=1$ if $m$ is a prime power and $Q=2$ otherwise, and $B_{1,\chi}$ is defined by $\tfrac1{f_{\chi}} \sum_{a=1}^{f_{\chi}} a\cdot\chi(a)$.%, and $B_{1,\chi}$ is defined by $\tfrac1f \sum_{a=1}^f a\cdot\chi(a)$ for any odd primitive character $\chi$ modulo $m$ of conductor~$f$.

Surprisingly, we observed that the determinant computation is very competitive, especially when the number of coprime factors of $m$ is small. Some comparative timings are provided in \cref{Table_hminus}.
\begin{table}[!htb]
   \centering
   \begin{tabular}{LLCRR}
     \hline

     \hline
     \multirow{2.5}{*}{$m$} & \multirow{2.5}{*}{$q_1\dotsc q_t$} & \multirow{2.5}{*}{$\varphi(m)$} & \multicolumn{2}{c}{Time $h_m^-$ (s)}\\[.3em]
     \cline{4-5}\\[-1.4\medskipamount]
     &&&\text{Analytic}  & \text{$[\A_m:\St m]$}\\
     \hline

     \hline
     1139 &17\cdot67                & 1056 & 12.6 & 8.1 \\
     1495 &5\cdot13\cdot23          & 1056 & 7.6  & 7.9 \\
     4140 &2^2\cdot3^2\cdot5\cdot23 & 1056 & 4.8  & 8.5 \\
     \hline
     2283 &3\cdot761	            & 1520 & 25.1 & 21.8 \\
     2865 &3\cdot5\cdot191	    & 1520 & 16.3 & 21.0 \\
     %3056 &2^4\cdot191              & 1520 & 18.1 & 18.1 \\
     \hline
     1951 &1951	                    & 1950 & 78.8 & 60.3 \\
     \hline
     2171 & 13\cdot167              & 1992 & 57.6 & 35.6 \\
     2495 & 5\cdot499               & 1992 & 53.8 & 41.7 \\
     %3507 & 3\cdot7\cdot167         & 1992 & 31.7 & 41.1 \\
     6012 & 2^2\cdot3^2\cdot167     & 1992 & 28.3 & 40.2 \\
     \hline
     
     \hline
   \end{tabular}
   \caption{Comparative timings for computing the relative class number~$h_m^-$ using resp. the analytic formula \cref{analytic_hminus} and the index formula for $[\A_m:\St m]$ in \cref{Relative_hK},  for a few representative examples.}
   \label{Table_hminus}
\end{table}

% A note on Jacobi sums generators
Finally, we verified that relations $\LL^{\alpha_m(b)} = J(a_1,a_2,\LL)\cdot\Z[\zeta_m]$ hold true in small dimensions (up to $\varphi(m)=80$). We note that computing explicitly such generators using the Jacobi sum formalism is very easy for any $m$. For instance, taking $m=2003$ and $\ell = 48073 \equiv 1\pmod m$, the computation of all $\varphi(m)/2$ generators corresponding to $\LL^{\alpha_m(b)}$, for all $b\in M'_m$ and some $\LL$ above $\ell$ takes under 15 minutes, \ie less than 1 second per generator.

By contrast, using suitable combinations of Gauss sums to obtain \eg generators for the $\LL^{(a-\si{a}{m})\cdot\theta_m(-1)}$ relations of \cite[\crefext{lem}{6.9}]{Was97} imposes to work in~$\Q\bigl[\zeta_{m\ell}\bigr]$. Even using all available algorithmic tricks, such as using sparse polynomials modulo $x^{m\ell}-1$, replacing divisions by~$g(b,\LL)\cdot g(-b,\LL) = \pm \norm{\LL}$ \cite[\crefext{lem}{6.1(b)}]{Was97} and profitting from \cref{Galois_action2}, this is arguably intractable in the above case when $m\ell = 96\,290\,219$, and still takes over 39 seconds per generator when restricting to the first split prime $\ell=4007$.

\begingroup
\small
\bibliographystyle{bibstyle}
\bibliography{Short_basis}
\endgroup

\end{document}